\documentclass[11pt]{amsart}

\usepackage{amssymb}
\usepackage{amsfonts}
\usepackage{amsmath}
\newtheorem{theorem}{Theorem} 
\newtheorem{lemma}[theorem]{Lemma}
\newtheorem{proposition}[theorem]{Proposition}
\newtheorem{corollary}[theorem]{Corollary}

\newtheorem{ques}{Question}

\newtheorem{preexample}[theorem]{Example} 

\makeatletter
\newcommand{\mysubsection}%
{\@startsection{subsection}{2}{\z@}{-3.25ex plus -1ex minus -.2ex}{-1ex}{\normalsize\sc}}
\makeatother

\renewcommand{\phi}{\varphi}
\renewcommand{\epsilon}{\varepsilon}
\newcommand{\U}{\mathbb{U}}

\newcommand{\C}{\mathbb{C}}

\newcommand{\cphi}{C_\phi}

\newcommand{\cphistar}{\cphi^*}

\newcommand{\dw}{{\omega}}

\newcommand{\htwo}{H^2}

\newcommand{\N}{\mathbb{N}}

\renewcommand{\Re}{{\rm Re\,}}

\newcommand{\spec}{\text{spec}}

\newcommand{\vstrut}{\rule{0in}{0.16in}}

\newcommand{\Sp}{\text{Sp}}
\newcommand{\ESp}{\text{Sp_e}}

\begin{document}



\title{Normal weighted composition operators on the Hardy space $H^2(\U)$}
\author{Paul S. Bourdon}
\address{Department of Mathematics\\  Washington and Lee University, Lexington VA 24450}
\email{pbourdon@wlu.edu}

\author{Sivaram K.\ Narayan$^*$}\thanks{$^*$This author thanks Central Michigan University for its support during his Fall 2009 sabbatical leave and Washington and Lee University for the hospitality extended during his sabbatical visit. }
\address{Department of Mathematics\\ Central Michigan University, Mount Pleasant MI 48859}
\email{sivaram.narayan@cmich.edu}

\begin{abstract} Let $\phi$ be an analytic function on the open unit disc $\U$ such that $\phi(\U)\subseteq \U$, and let $\psi$ be an analytic function on $\U$ such that the weighted composition operator $W_{\psi, \phi}$ defined by $W_{\psi, \phi} f = \psi f \circ \phi$ is bounded on the Hardy space $H^2(\U)$.  We characterize those weighted composition operators on $H^2(\U)$ that are unitary, showing that in contrast to the unweighted case ($\psi \equiv 1$), every automorphism of $\U$ induces a unitary weighted composition operator.   A conjugation argument, using these unitary operators, allows us to describe all normal weighted composition operators on $H^2(\U)$ for which the inducing map $\phi$ fixes a point in $\U$.   This description shows both $\psi$ and $\phi$ must be linear fractional in order for $W_{\psi, \phi}$ to be normal (assuming $\phi$ fixes a point in $\U$).   In general, we show that if  $W_{\psi, \phi}$ is normal on $H^2(\U)$ and $\psi\not\equiv 0$, then $\phi$ must be either univalent on $\U$ or constant.     Descriptions of spectra are provided for the operator  $W_{\psi, \phi}: H^2(\U)\rightarrow H^2(\U)$ when it is unitary or when it is normal and $\phi$ fixes a point in $\U$. 
\end{abstract}

\maketitle 
 
 \section{Introduction}
  
   Let $H(\U)$ denote the collection of all holomorphic functions on the open unit disc $\U$, and let $H^2(\U)$ be the classical Hardy space of $\U$, a Hilbert space consisting of those functions in $H(\U)$ whose Maclaurin coefficients are square summable.  Throughout this paper, $\phi$ denotes an analytic selfmap of $\U$; i.e., an element of $H(\U)$ mapping $\U$ into $\U$.  A result due to Littlewood \cite{Lit} shows that  $\phi$ induces a bounded composition operator on $H^2(\U)$ defined by $C_\phi f = f\circ \phi$.   Here, we study weighted composition operators $W_{\psi, \phi}$ on $H^2(\U)$, which result from following composition with $\phi$ with multiplication by a weight function $\psi\in H(U)$; thus, $W_{\psi, \phi} f = \psi\, f \circ \phi$.   Such weighted composition operators are clearly bounded on $H^2(\U)$ when $\psi$ is bounded on $\U$, but the boundedness of $\psi$ on $\U$ is not necessary for $W_{\psi, \phi}$ to be bounded. For example, if $\phi(\U)$ has closure contained in $\U$, then $\psi \in H^2(\U)$ is sufficient to ensure $W_{\psi,\phi}: H^2(\U) \rightarrow H^2(\U)$ is bounded.  We call $W_{\psi,\phi}$ the weighted composition operator induced by $\phi$ with weight function $\psi$.   
   
    Normal composition operators on the Hardy space $H^2$ of the open unit disc are easily characterized.  In fact, merely the assumption that $C_\phi$ and $C_\phi^*$ commute when applied to the Hardy space functions $f(z) = 1$ and $f(z) = z$ reveals that if $C_\phi: H^2(\U)\rightarrow H^2(\U)$ is normal, then $\phi(z) = cz$ for some constant $c$, $|c| \le 1$ (see \cite[Exercise 8.1.2]{CMB}).   Conversely, if $\phi(z) = cz$, then $C_\phi$ is a diagonal operator relative to the orthonormal basis $(z^n)_{n=0}^\infty$ of $H^2(\U)$ and thus $C_\phi$ is normal.  
    
     For weighted composition operators, the normality question is much more interesting.  Indeed, this paper is inspired by Cowen and Ko's characterization of Hermitian weighted composition operators \cite{CK},  which reveals that with appropriate restrictions on $a_0$ and $a_1$, the pair  $\psi(z) = c/(1-\overline{a_0}z)$, $\phi(z) = a_0 + a_1z/(1-\overline{a_0}z)$ provide a weighted composition operator $W_{\psi, \phi}$ that is Hermitian whenever $c$ is real.    
    
    In this paper, we characterize those weighted composition operators on $H^2(\U)$ that are unitary and apply this characterization to describe all normal weighted composition  operators on $H^2(\U)$ that are induced by a selfmap $\phi$ of $\U$ which fixes a point in $\U$.   Our characterization of unitary weighted composition operators, which occupies Section~\ref{USect},  shows that every automorphism $\phi$ of $\U$  has a companion weight function $\psi$ such that $W_{\psi, \phi}$ is unitary on $H^2(\U)$.   In contrast, only the rotation maps $\phi(z) = \zeta z$, $|\zeta| = 1$, induce unitary composition operators $C_\phi: H^2(\U)\rightarrow H^2(\U)$.     Our description of all normal weighted composition operators on $H^2(\U)$ induced by selfmaps $\phi$ of $\U$ fixing a point in $\U$, which appears in Section~\ref{SNWCO}, reveals that for such operators both $\psi$ and $\phi$ are linear fractional.  In the ``Preliminaries'' section below, we establish that if  $W_{\psi,\phi}$  is  normal on $H^2(\U)$ and $\psi\not\equiv 0$, then either  $\phi$ is univalent on $\U$ or constant and $\psi(z)$ is nonzero at every $z\in \U$.   At the conclusions of Sections~\ref{USect} and \ref{SNWCO}, we present spectral characterizations of the unitary and normal operators identified in the section.
    
     In the final section of the paper, we  point out that the weight functions $\psi$ that give rise  to  weighted composition operators $W_{\psi, \phi}$ that are unitary, Hermitian, or normal (with inducing map fixing a point of $\U$) all have the same form, and then determine when such $\psi$ can pair with a linear fractional $\phi$ to produce a normal operator $W_{\psi, \phi}$.    

   \section{Preliminaries}\label{PSect}
    
    Here we collect some background information necessary to our work and then present some simple necessary conditions for $W_{\psi, \phi}:H^2(\U)\rightarrow H^2(\U)$  to be normal.
    
      \subsection{Reproducing kernels for $H^2(\U)$}   The Hardy space $H^2(\U)$ is a Hilbert space with inner product 
$$
\langle f, g\rangle = \sum_{n=0}^\infty \hat{f}(n)\overline{\hat{g}(n)},
$$
 where $(\hat{f}(n))$ and $(\hat{g}(n))$ are the sequences of Maclaurin coefficients for $f$ and $g$ respectively.  The norm of $f\in H^2(\U)$ is given by $\left(\sum_{n=0}^\infty |\hat{f}(n)|^2\right)^{1/2}$ or, alternatively, by
\begin{equation}\label{BIN}
\|f\|_{H^2(\U)}^2 = \frac{1}{2\pi} \int_0^{2\pi} |f(e^{it})|^2\, dt,
\end{equation}
where $f(e^{it})$ represents the radial limit of $f$ at $e^{it}$, which exists for a.e. $t\in [0, 2\pi)$.  

Let $h\in H^{\infty}(\U)$, the Banach algebra of bounded analytic functions on $\U$ with $\|h\|_\infty = \sup\{|h(z)|: z\in \U\}$.   The integral representation of the $H^2(\U)$ norm makes it clear that  $\|hf\|_{H^2(\U)} \le \|h\|_\infty \|f\|_{H^2(\U)}$ for each $f\in H^2(\U)$.  Thus, the {\it multiplication operator} $M_h: H^2(\U) \rightarrow H^2(\U)$, defined by $M_h f = hf$, is bounded and linear on $H^2(\U)$.  When $\psi \in H^\infty(\U)$, note that we may write $W_{\psi, \phi} = M_{\psi} C_\phi$ and $\|W_{\psi, \phi}  \| \le \|M_\psi\| \|C_\phi\| = \|\psi\|_\infty \|C_\phi\|$.

For each $\beta\in \U$, let $K_\beta = 1/(1-\bar{\beta}z)$.  Then $K_\beta \in H^2(\U)$, and it is easy to see that  for each function $f\in H^2(\U)$, 
$$
\langle f, K_\beta\rangle = f(\beta);
$$
thus, $K_\beta$ is the {\it reproducing kernel} at $\beta$ for $H^2(\U)$.  We make extensive use of the following simple and well-known lemma.

\begin{lemma}\label{WRKL} Suppose that $W_{\psi, \phi}:H^2(\U) \rightarrow H^2(\U)$ is bounded and $\beta\in \U$.  Then
$$
W_{\psi, \phi}^*  K_\beta = \overline{\psi(\beta)} K_{\phi(\beta)}.
$$
\end{lemma}
\begin{proof}  Let $f$  be an arbitrary function in $H^2(\U)$. We have
$$
\langle f, W_{\psi,\phi}^* K_\beta \rangle  =  \langle \psi f\circ \phi, K_\beta\rangle \\
=  \psi(\beta)f(\phi(\beta)) \\
 =  \langle f, \overline{\psi(\beta)} K_{\phi(\beta)}\rangle
$$
and the lemma follows.
\end{proof}

\subsection{The Denjoy-Wolff Point $\dw$.}   Let $\phi$ be an analytic selfmap of $\U$.  For $n$ a nonnegative integer, let $\phi^{[n]}$ denote the $n$-th iterate of $\phi$ so that, e.g.,  $\phi^{[0]}$ is the identity function on $\U$ and $\phi^{[2]} = \phi \circ \phi$.     If $\phi$ is not an elliptic automorphism of $\U$, then there is a  (unique)  point $\dw$ in the closure $\U^-$ of $\U$ such  that
$$
\dw = \lim_{n\rightarrow \infty} \phi^{[n]}(z)
$$
for each $z\in \U$.   The point $\dw$, called the {\it Denjoy-Wolff point} of $\phi$, is also
characterized as follows: if $|\dw| < 1$, then $\phi(\dw) = \dw$ and
$|\phi'(\dw)| < 1$; if
$\dw\in \partial \U$, then $\phi(\dw) =\dw$ and $0 < \phi'(\dw) \le 1$.  If
$|\dw| = 1$, then $\phi(\dw)$ represents the angular (non-tangential) limit of $\phi$
at
$\dw$ and $\phi'(\dw)$ represents the angular derivative of $\phi$ at $\dw$.    The location of the Denjoy-Wolff point and the behavior of iterate sequences $(\phi^{[n]}(z))$ as they approach the Denjoy-Wolff point strongly influence properties of the weighted composition operator $W_{\psi, \phi}$.  To see how spectral behavior varies with location of $\dw$, see, e.g, \cite[Section 5]{BrdS}.

     \subsection{Cowen's formula for the adjoint of a linear fractional composition operator}  
     
   As we have indicated, we will show that when $W_{\psi, \phi}$ is unitary or when it is normal and $\phi$ fixes a point in $\U$, the inducing map $\phi$ must be linear fractional and the weight function $\psi$ must be both linear fractional and bounded on $\U$.   Thus $W_{\psi, \phi} ^*= (M_\psi C_\phi)^* = C_\phi^* M_\psi^*$, and hence Cowen's formula for the adjoint of a linear fractional composition operator becomes quite important to our work.  
   
    When $\phi(z) = (az + b)/(cz + d)$ is a nonconstant linear fractional selfmapping of $\U$, Cowen \cite{Cow} establishes 
   \begin{equation}\label{CF}
    C_\phi^* = M_g C_\sigma M_h^*
\end{equation}
    where the {\it Cowen auxiliary functions}  $g$, $\sigma$ and $h$ are defined as follows:
    $$
    g(z) = \frac{1}{\vstrut-\bar{b}z + \bar{d}}\, , \quad \sigma(z) = \frac{\bar{a}z-\bar{c}}{\vstrut -\bar{b}z + \bar{d}}\, ,\  \text{and}\quad  h(z) = cz + d.
    $$
\subsection{From the unit disk to the right halfplane and back}  The function $T(z) = (1+z)/(1-z)$ maps the unit disk $\U$ univalently onto the right halfplane $\Pi$.   If $\phi$ is selfmap of $\U$ with Denjoy-Wolff point $1$, then $\Phi: = T\circ \phi\circ T^{-1}$ is a selfmap of $\Pi$ with Denjoy-Wolff point infinity (meaning $\Phi^{[n]}(w) \rightarrow \infty$ as $n\rightarrow \infty$ for every $w\in \Pi$).  

In the sequel, we will use generic formulas for linear fractional selfmappings of $\U$ having Denjoy-Wolff point $1$.  These formulas are easily derived using the correspondence between $\U$ and $\Pi$ provided by $T$.  

Suppose that $\phi$ is a linear fractional selfmap of $\U$ with Denjoy-Wolff point $1$ such that $\phi(1) = 1$ and $\phi'(1) = 1$.  (Thus $\phi$ is of  {\it parabolic type}.)  It is easy to see that $\Phi(w) = (T\circ\phi\circ T^{-1})(w) = w + t$, where $\Re(t) \ge 0$ and is nonzero.  (Note $\Re(t) = 0$ if and only if $\phi$ is an automorphism.)  Hence, 
\begin{equation}\label{PF}
\phi(z) = T^{-1}(T(z) + t) = \frac{(2-t)z + t}{-tz + (2+t)}.
\end{equation}

Suppose that $\phi$ is a linear fractional selfmap of $\U$ such that $\phi(1) = 1$ and $\phi'(1)  = b < 1$.  (Thus $\phi$ is of {\it hyperbolic type}.) It is easy to see that $\Phi(w) = (T\circ\phi\circ T^{-1})(w) = rw + t$, where $r = 1/b$ and $\Re(t) \ge 0$ (and $\Re(t) = 0$ if and only if $\phi$ is an automorphism).   Hence, 
\begin{equation}\label{HF}
\phi(z) = T^{-1}(T(z) + t) = \frac{(1+r-t)z + r+t -1}{(r-t-1)z+1 +r+t}.
\end{equation}

\subsection{Two necessary conditions for normality of $W_{\psi, \phi}$}

 \begin{lemma}\label{MNZ} If  $W_{\psi, \phi}$ is normal then either $\psi\equiv 0$ or $\psi$ never vanishes on $\U.$
 \end{lemma}
  \begin{proof}   Suppose $W_{\psi, \phi}$ is normal and $\psi(\beta)=0$ for some $\beta$ in $\U.$ Then by Lemma~\ref{WRKL},  $W_{\psi,\phi}^*K_\beta = \overline {\psi(\beta)} K_{\phi(\beta)} \equiv 0.$ Since $W_{\psi, \phi}$ is normal,  $\|W_{\psi,\phi} K_\beta\|= \|W_{\psi, \phi}^*K_\beta\|.$  Therefore, $\|W_{\psi,\phi} K_{\beta}\| =  0$ and thus $\frac{\psi(z)}{1-\overline{\beta} \phi(z)} = 0$ for every $z$ in $\U$, which implies $\psi\equiv 0$.  Thus, if $W_{\psi, \phi}$ is normal either $\psi\equiv 0$ or $\psi$ is nonzero at each point in $\U$.
\end{proof}

\begin{proposition}\label{WUPU} Suppose $W_{\psi, \phi}$ is normal.  If $\phi$ is not a constant function and $\psi$ is not the zero function then $\phi$ is univalent.
\end{proposition}
\begin{proof} Suppose $\phi$ is not univalent on $\U$ and is nonconstant. Then there exist points $a$ and $b$ in $\U$ such that $a\neq b$ and $\phi(a) = \phi(b)$. Since $\psi \not\equiv 0$, from Lemma~\ref{MNZ}, we conclude that $\psi(a) \neq 0$ and $\psi(b) \neq 0$. Let 
$$
g = \frac{K_a}{\vstrut \overline{\psi(a)}} - \frac{K_b}{\vstrut \overline{\psi(b)}}
$$
and observe $g$ is a nonzero function in $H^2(\U)$.  We have $W_{\psi, \phi}^*g = 0$ by Lemma~\ref{WRKL}. Since  $W_{\psi, \phi}$ is normal, $\|W_{\psi, \phi} g\| = \|W_{\psi, \phi}^* g\|=0$.  But $\|W_{\psi, \phi} g\| = 0$ implies $g\circ \phi$ is the zero function.  Since $\phi$ is nonconstant, $g$ must vanish on a nonempty open subset of $\U$ and hence  $g$ must be the zero function, a  contradiction.  Hence $\phi$ is univalent.
\end{proof}

\section{Unitary weighted composition operators}\label{USect}

        Suppose that $\phi$ is the constant function $z\mapsto \beta$ for $\beta\in \U$.  Then  $W_{\psi, \phi}$ annihilates any function in $H^2(\U)$ that vanishes at $\beta$ and hence $W_{\psi, \phi}$ cannot be norm-preserving and thus cannot be unitary.  Hence if $W_{\psi,\phi}$ is unitary, we may apply Proposition~\ref{WUPU} to conclude that $\phi$ must be univalent on $\U$.   In fact, $\phi$ must be an automorphism of $\U$.
        
        \begin{proposition}\label{UPU} Suppose that $W_{\psi, \phi}:H^2(\U)\rightarrow H^2(\U)$ is unitary.  Then $\phi$ must be an automorphism of $\U$.
        \end{proposition}
        \begin{proof}   Because unitary operators are norm preserving,  we see 
        $$
        1 = \|W_{\psi,\phi} 1\| = \|\psi\|.
        $$
        In addition, if $f(z) = z$, then $\|f\| = 1$ and 
        $$
        1 = \|W_{\psi, \phi} f\| = \| \psi \phi\|.
        $$
Since $|\phi(e^{it})| \le 1$ for a.e.\ $ t\in [0, 2\pi)$ and  both $\|\psi\|$ and $\|\psi \phi\|$ are $1$, the integral representation \ref{BIN} of the norm on $H^2(\U)$ shows that $|\phi(e^{it})| = 1$ a.e. on $\U$ so that $\phi$ is an inner function.  However, as we have indicated, thanks to Proposition~\ref{WUPU}, we know $\phi$ is univalent.   It is not difficult to see that a univalent inner function must be an automorphism of $\U$ (see, e.g., \cite[Corollary 3.28, p.\ 152]{CMB} ), and the proposition follows.
\end{proof}

Now we identify what form the multiplier $\psi$ must take in order that $W_{\psi, \phi}$  be unitary.  Because any automorphism $\phi$ of $\U$ must take the value $0$ at some $\beta$ in $\U$, the following proposition shows that if $W_{\psi, \phi}$ is unitary, then $\psi = ck_{\phi^{-1}(0)}$, where $|c| = 1$ and $k_{\phi^{-1}(0)}$ is the normalized reproducing kernel  $K_{\phi^{-1}(0)}/\|K_{\phi^{-1}(0)}\|$.

\begin{proposition}\label{UMF} Suppose the inducing map $\phi$ satisfies $ \phi (\beta) = 0$ for some $\beta\in \U$.  If  $W_{\psi, \phi}$ is unitary,  then $$\psi = c \frac {K_\beta}{\|K_{\beta}\|}$$ where $|c| =1.$
 \end{proposition}
 
 \begin{proof}  Suppose that $W_{\psi, \phi}$ is unitary;  then $W_{\psi,\phi}W_{\psi,\phi}^*K_\beta = K_\beta$.  Applying Lemma~\ref{WRKL},   we can rewrite the preceding equation as 
 $$
 W_{\psi,\phi} \overline{\psi(\beta)} K_0 = K_\beta.
 $$
   Because $K_0\equiv 1$, we obtain  $\psi \overline{\psi(\beta)} = K_\beta$, so that 
  $$
  \psi = \frac{K_\beta}{\overline{\psi(\beta)}}.
  $$
        Evaluating both sides of the preceding equation at $\beta$ yields  $|\psi(\beta)|^2 = K_\beta(\beta) = \|K_\beta\|^2$ and the proposition follows.
        \end{proof}

\begin{theorem}\label{UWCO} The weighted composition operator $W_{\psi, \phi}$ is unitary on $H^2(\U)$ if and only if  $\phi$ is an automorphism of $\U$ and $\psi = c K_\beta/\|K_\beta\|$ where $\phi(\beta) = 0$ and $|c| =1 .$
\end{theorem}

\begin{proof}  If $W_{\psi, \phi}$ is unitary then $\phi$ is an automorphism of $\U$ by Proposition~\ref{UPU} and $\psi$ has the form claimed by Proposition~\ref{UMF}.    

\par To prove sufficiency, let 
$$
\phi(z) =  \eta\frac {\beta - z}{1 - \overline{\beta}z}, \quad  \text{and} \quad  \psi = c \frac {K_\beta}{\|K_\beta\|}, 
$$
where $|\eta| =1$ and  $|c| =1$. The Cowen auxiliary functions of $\phi$ are
$$
\sigma(z) = \frac{\beta-\bar{\eta} z} {\vstrut 1 - \overline{\eta \beta }z} = {\phi}^{-1}(z),  \
g(z) = \frac{1}{1 - \overline {\eta \beta}z}, \  \text{and}\  h(z) = 1 - \overline{\beta}z.
$$
Observe that $g\circ \phi = \|K_\beta\|^2/K_\beta$ and $\psi\circ \sigma = c\|K_\beta\|/g$.

Since $M_{h}^* M_{\psi}^* = \overline{c}/\|K_b\|$ is a constant, $ \sigma = {\phi}^{-1}$, and $W_{\psi, \phi} = M_\psi C_\phi$ (because $\psi$ is bounded), we see that
$$
W_{\psi, \phi} W_{\psi, \phi}^*= M_{\psi} C_{\phi} M_{g} C_{\sigma} M_{h}^* M_{\psi}^* =  \frac {\overline c}{\|K_\beta\|} \cdot \psi\cdot  g\circ \phi \cdot C_{\sigma \circ \phi} = I
$$
and
$$
W_{\psi, \phi}^* W_{\psi, \phi} = M_{g} C_{\sigma} M_{h}^* M_{\psi}^* M_{\psi} C_{\phi} = \frac {\overline c}{\|K_\beta\|}\cdot  g\cdot \psi \circ\sigma \cdot C_{\phi\circ\sigma} = I.
$$
This completes the proof of the theorem.
\end{proof}

Unitary weighted composition operators $W_{\psi, \phi}$ divide naturally into three classes based on whether the automorphism $\phi$ is (i) elliptic, (ii) hyperbolic, or (iii) parabolic.    This observation is the  basis for our spectral characterizations of these operators.   
\begin{theorem}\label{SCT} Suppose that $W_{\psi,\phi}$ is unitary.  
\begin{itemize}
\item[(a)] If $\phi$ is elliptic, fixing $p \in \U$, then $|\phi'(p)|=1$ and the spectrum of $W_{\psi, \phi}$ is the closure of  $\{\psi(p)\phi'(p)^n: n=0,1, 2, \ldots\}$.  
\item[(b)] If $\phi$ is hyperbolic or parabolic then the spectrum of $W_{\psi, \phi}$ is the unit circle.
\end{itemize}
\end{theorem}
\begin{proof}   Suppose that $W_{\psi, \phi}$ is unitary so that $ \psi = c  K_\beta/\|K_\beta\|$, where $\phi(\beta) = 0$ and $\phi$ is an automorphism.   We know the spectrum of $W_{\psi,\phi}$ must be contained in the unit circle.   

In case (a),  $W_{\psi, \phi}$ is a normal operator whose inducing map $\phi$ fixes a point $p$  in $\U$.  That the spectrum of $W_{\psi,\phi} = \{\psi(p)\phi'(p)^n: n=0,1, 2, \ldots\}$  follows immediately from Proposition~\ref{SCNIFP} of Section~\ref{SNWCO} below, which  characterizes the spectrum of any normal weighted composition operator on $H^2(\U)$ induced by a selfmap of $\U$ fixing a point in $\U$.  We wish to make clear here the connection between the form of our unitary weighted composition operators in this elliptic case and the form~(\ref{IFPN}) of the normal  weighted compositions described in the next section.  By doing so, we will confirm that $|\psi(p)\phi'(p)| = 1$, which must be true since the spectrum of $W_{\psi, \phi}$ is a subset of the unit circle.    We are assuming $\phi(p) = p$. We can conjugate $\phi$ by the self-inverse automorphism $\alpha_p(z) = (p-z)/(1-\bar{p}z)$ of $\U$, obtaining that $\alpha_p \circ\phi\circ \alpha_p = \zeta z$ for some unimodular constant $\zeta$; in fact, $\zeta = \phi'(p)$.   Equivalently, $\phi = \alpha_p\circ(\zeta\alpha_p)$.  Since $\phi(\beta) = 0$, we see $\beta = \alpha_p(\bar{\zeta}p)$.    Now observe that $\psi = c K_\beta/\|K_\beta\|$ is also given by
$$ 
\psi(z) = c \frac{1-|p|^2\zeta}{\vstrut |1-|p|^2 \bar{\zeta}|}\frac{K_p}{K_p\circ \phi}.
$$
Note $\psi(p) = c (1-|p|^2\zeta)/|1-|p|^2\bar{\zeta}|$ is a unimodular constant.    Thus in case (a) we have 
$$
\psi = \psi(p)\frac{K_p}{K_p\circ \phi} \quad  \text{and} \quad \phi = \alpha_p\circ (\zeta\alpha_p).
$$   

 Turning to case (b), let's initially assume that $\phi$ is of hyperbolic type, which means that $\phi$ has Denjoy-Wolff point $\dw \in \partial \U$ and $\phi'(\dw) <1$.   We employ the idea of the proof  \cite[Theorem 4.3] {Cow2}  to see  that $W_{\psi, \phi}$ similar to $\zeta W_{\psi, \phi}$ for every unimodular constant $\zeta$.         Thus, the spectrum of $W_{\psi, \phi}$, which we already knew to be a subset of $\partial \U$, must be all of $\partial \U$.   To see that  $W_{\psi, \phi}$ is similar to $\zeta W_{\psi, \phi}$, take $\zeta = e^{i\theta}$.  To establish similarity, one can use the function $f$ identified in the proof of \cite[Theorem 4.3]{Cow2}: this function satisfies $f, 1/f \in H^\infty(\U)$,  and $C_\phi f = e^{i\theta} f$.    It is easy to see that $M_f^{-1} W_{\psi, \phi} M_f = \zeta W_{\psi,\phi}$.
 
   Now, we assume that $\phi$ is a parabolic automorphism and $W_{\psi, \phi}$ is normal.  We know by Theorem~\ref{UWCO} that $\psi = cK_\beta/\|K_\beta\|$, where $|c| = 1$ and  $\phi(\beta) = 0$.   Because $\phi$ is a parabolic automorphism,  $\phi$ has Denjoy-Wolff point $\dw \in \partial \U$ and $\phi'(\dw) =1$.   Without loss of generality,  we may assume that $\dw = 1$. (If $\dw\ne 1$, let $U = C_{\dw z}$ so that $U$ is unitary and note $UW_{\psi,\phi} U^* = W_{\tilde{\psi}, \tilde{\phi}}$, where $\tilde{\phi}(z) = \bar{\dw}\phi(\dw z)$ is a parabolic automorphism with fixed point $1$ and $\tilde{\psi} = K_{\tilde{\beta}}/\|K_{\tilde{\beta}}\|$ where $\tilde{\phi}(\tilde{\beta}) = 0$.) Since $\phi$ is a parabolic automorphism with $\phi(1) = 1$,  We may assume that $\phi$ has the form given by (\ref{PF}) with $\Re(t) = 0$ and $t\ne 0$:
\begin{equation}\label{FOPP}
   \phi(z) =  \frac{(2-t)z + t}{-tz + (2+t)}.
\end{equation}
 We know the spectrum of $W_{\psi,\phi}$ is contained in the unit circle.  That it equals the unit circle is a consequence of Theorem 19 of \cite{CK}, or,  more accurately, it follows from an argument just like that of  \cite[Theorem 19]{CK}.   For completeness we include a proof.  Let $\lambda < 0$ be arbitrary.  We show that there is a unimodular constant $\nu$ such that
 $$
 W_{\psi, \phi} - \nu e^{\lambda t}I 
 $$
 is not bounded below on $H^2(\U)$, which, since $t$ nonzero and pure imaginary and $\lambda < 0$ is arbitrary, will complete the proof.    
 
 We are assuming $\phi$ is given by (\ref{FOPP}), so that $\phi(\beta) = 0$ where $\beta = t/(t-2)$.  Its companion weight function has the form
 $$
 \psi = c \frac{K_\beta}{\|K_\beta\|} =  \frac{2c}{\sqrt{4+ |t|^2}} \frac{2+t}{2+t-tz},
 $$
 where we have used $\bar{t} = - t$ since $t$ is pure imaginary.    The proof is based on the observation that $G(z) = \frac{1}{1-z}\exp\left(\lambda \frac{1+z}{1-z}\right)$ satisfies
 $$
 W_{\psi, \phi} G = \frac{c(2+t)}{\sqrt{4+|t|^2}} e^{\lambda t}G.
 $$
 Let $\nu = c(2+t)(4+|t|^2)^{-1/2}$ so that $\nu$ is a unimodular constant and $W_{\psi, \phi} G =\nu e^{\lambda t} G$.  It appears that the point spectrum of $W_{\psi,\phi}$ contains the unimodular constant $\nu e^{\lambda t} $ for every $\lambda < 0$.  However, $G$ is not in $H^2$.  The idea of the proof of \cite[Theorem 19]{CK} is to modify $G$ as follows 
 $$
 G_s(z) = \frac{1}{s-z}\exp\left(\lambda \frac{1+z}{1-z}\right),   $$
so that  $G_s\in H^2(\U)$ for every real number $s$ such that $s > 1$, and then  show
\begin{equation}\label{BBC}
\left\| W_{\psi, \phi} \frac{G_s}{\|G_s\|} -  \nu e^{\lambda t}\frac{G_s}{\|G_s\|}\right\| \rightarrow 0 \quad {\rm as} \quad  s\rightarrow 1^+.
\end{equation}
Because $S_\lambda(z):=\exp\left(\lambda \frac{1+z}{1-z}\right)$ is an inner function,  $\|G_s\|^2 = \|1/(s-z)\|^2 = s^{-2}\| K_{1/s}\|^2 = 1/(s^2 - 1)$. A computation shows
$$
W_{\psi, \phi} G_s =  e^{\lambda t}  2\nu \frac{1}{(t-st-2)z +  2s + st - t} \exp\left(\lambda \frac{1+z}{1-z}\right),
$$
Thus, since $\nu$ is unimodular, $t$ is pure imaginary, and $S_\lambda$ is inner, we see the norm on the left of (\ref{BBC}) simplifies to 
$$
\sqrt{s^2-1} \left\| \frac{2}{(t-st-2)z  +2s + st - t} - \frac{1}{s-z}\right\| = \sqrt{s^2 -1} \left \| \frac{2}{2s+st - t}K_{\frac{2-st+t}{2s-st+t}} - \frac{1}{s}K_{1/s}\right\|,
$$
where we used $\bar{t} = -t$  in obtaining the reproducing kernel representation of the fraction on the left within the norm.
Expanding 
$$
\left\langle \frac{2}{2s+st - t}K_{\frac{2-st+t}{2s-st+t}} - \frac{1}{s}K_{1/s}, \frac{2}{2s+st - t}K_{\frac{2-st+t}{2s-st+t}} - \frac{1}{s}K_{1/s}\right\rangle
$$
yields  (since $t$ is pure imaginary)
\begin{equation*}
\begin{split}
\frac{4}{|2s + st -t|^2-|2-st +t|^2} + \frac{1}{s^2 - 1} - 4\frac{1}{(s-1)}\Re\left( \frac{1}{(2(s+1)-st + t)}\right)  \\= \frac{2}{s^2-1} - \frac{8(s+1)}{s-1}\frac{1}{4(s+1)^2 + |t|^2(1-s)^2}.
\end{split}
\end{equation*}
Multiplying the quantity on the right of the preceding equation by $s^2-1$ and taking the limit as $s\rightarrow 1^+$ yields $2-2 = 0$ and we have established (\ref{BBC}), and hence completed our proof that the spectrum of $W_{\psi, \phi}$ is the unit circle.\end{proof}

\section{Normal $W_{\psi, \phi}$  with  $\phi$ fixing a point of $\U$}\label{SNWCO}

If $W_{\psi, \phi}:H^2(\U)\rightarrow H^2(\U)$ is normal and $\phi(p) = p$ for some $p\in \U$, then it is easy to prove that the weight function $\psi$ has a simple, linear fractional form.

\begin{proposition}\label{FOM}  Suppose that $\phi$ has a fixed point  $p\in \U$ and $W_{\psi, \phi}$ is normal.  Then 
$$
\psi =  \psi(p) \frac {K_{p} }{K_{p} \circ \phi}.
$$
\end{proposition}
\begin{proof}   We have 
$$
W_{\psi ,\phi}^* K_p = \overline{\psi(p)} K_p
$$
so that $K_p$ is an eigenvector for $W_{\psi, \phi}^*$ with corresponding eigenvalue $\overline{\psi(p)}$.  Since $W_{\psi,\phi}$ is normal $K_p$ is an eigenvector for $W_{\psi, \phi}$ with corresponding eigenvalue $\psi(p)$ and we must have
\begin{equation}\label{NE}
\psi(p) K_p = W_{\psi, \phi} K_p = \psi \, K_p\circ \phi,
\end{equation}
from which the proposition follows.
\end{proof}

Because $K_0 \equiv 1$, the preceding proposition immediately yields:

\begin{corollary}\label{NWCOWO}  Suppose that  $\phi(0) = 0$.  Then $W_{\psi,\phi}$ is normal if and only if $\psi$ is constant and $C_\phi$ is normal.
\end{corollary}

\begin{theorem}\label{CNWCO} Suppose that $\phi$ has a fixed point  $p\in \U$.  Then $W_{\psi, \phi}$ is normal if and only if 
\begin{equation}\label{IFPN}
\psi = \gamma \frac {K_{p} }{K_{p} \circ \phi} \quad     and   \quad   \phi = {\alpha}_p \circ ({\delta {\alpha_p}})
\end{equation}
where $\alpha_{p}(z) = (p - z)/(1 - \overline{p} z)$  and  $\delta$  is a constant satisfying $|\delta| \le 1$ while $\gamma$ is a constant giving the value of $\psi$ at $p$ .
\end{theorem}
\begin{proof}  Suppose that $W_{\psi, \phi}$ is normal. Let $ \psi_{p} =  \frac {K_{p}}{\|K_{p}\|}$.  Since $\alpha_{p}$ is an automorphism taking $p$ to $0$ we know from Theorem~\ref{UWCO} that $W_{\psi_{p}, \alpha_{p}}$ is unitary. It follows that the operator
$$
L := (W_{\psi_{p}, \alpha_{p}}^*) (W_{\psi, \phi}) (W_{\psi_{p}, \alpha_{p}})
$$
is normal.
  Let $\sigma$,  $g$ and $h$ be the Cowen auxiliary functions for $\alpha_p$ so that $\sigma = \alpha_p^{-1} = \alpha_p$ and
  $h(z) = 1-\bar{p}z$,  $g(z) = 1/(1-\bar{p}z)$.  Observe 
  \begin{equation}\label{UPN}
  C_{\alpha_p}^* = M_g C_{\alpha_p} M_h^* \quad {\rm and} \quad   M_h^*M_{\psi_p}^* = (M_{\psi_p}M_h)^* = M_{1/\|K_p\|}.
  \end{equation}
     Setting $\mu = 1/\|K_p\|$ and using (\ref{UPN}), we find 
  $$
  L = C_{\alpha_p}^* M_{\psi_p}^* M_\psi C_\phi M_{\psi_p}C_{\alpha_p}\\
      =  \mu M_g M_{\psi\circ \alpha_p} M_{\psi_p\circ \phi\circ\alpha_p} C_{\alpha_p \circ\phi\circ\alpha_p} \\
     =  M_q C_{\alpha_p\circ \phi\circ \alpha_p},
$$
where $q\in H^\infty(\U)$ is given by $q = \mu g\cdot \psi\circ \alpha_p \cdot \psi_p\circ \phi\circ\alpha_p$.  
Since $L = W_{q, \alpha_p\circ \phi\circ \alpha_p}$ is normal and its inducing map $\alpha_p\circ \phi\circ \alpha_p$ fixes $0$, we may apply Corollary~\ref{NWCOWO} to conclude that $q$ must be a constant map and $C_{\alpha_p\circ \phi\circ \alpha_p}$ must be normal.  As we noted in the Introduction, because $C_{\alpha_p\circ \phi\circ \alpha_p}$ is normal on $H^2(\U)$, there must be a constant $\delta$ with $|\delta|\le 1$ such that $\alpha_p\circ \phi\circ \alpha_p = \delta z$ and we see 
$\phi = \alpha_p\circ(\delta \alpha_p)$, as desired. We can confirm that $q = \mu g\cdot \psi\circ \alpha_p \cdot \psi_p\circ \phi\circ\alpha_p$  is constant and determine its value by a straightforward computation.  We know from Proposition~\ref{FOM} that $\psi = \psi(p)K_p/\|K_p\|$;  substituting this formula for $\psi$ into  $\mu g\cdot \psi\circ \alpha_p \cdot \psi_p\circ \phi\circ\alpha_p$ and simplifying yields $q\equiv \psi(p)$.
  
     We now must show that if $\psi$ and $\phi$ have the forms specified in (\ref{IFPN}), then $W_{\psi, \phi}$ is normal.       We have
\begin{equation}\label{AE}
\phi(z) = \alpha_p(\delta\alpha_p(z)) =  \frac{p(1-\delta)+(\delta - |p|^2)z}{1-|p|^2\delta + \bar{p}(\delta-1)z}.
\end{equation}
Without loss of generality, we take $\gamma =1$ so that
\begin{equation}\label{phiF}
\psi(z) = \frac{K_p(z)}{K_p(\phi(z))} = \frac{1-|p|^2}{1-|p|^2\delta - \bar{p}(1-\delta)z}.
\end{equation}
We show that 
\begin{equation}\label{FCRY}
W_{\psi_p, \alpha_p} C_{\delta z} W_{\psi_p, \alpha_p}^* = W_{\psi, \phi}  
\end{equation}
 so that $W_{\psi, \phi}$ is unitarily equivalent to the normal operator $C_{\delta z}$ and hence is normal.
Using (\ref{UPN}), we see
$$
W_{\psi_p, \alpha_p}^* = \frac{1}{\|K_p\|} M_g C_{\alpha_p},
$$
where $g(z) = 1/(1-\bar{p}z)$.  Thus
\begin{equation}\label{FCRY2}
W_{\psi_p, \alpha_p} C_{\delta z} W_{\psi_p, \alpha_p}^* = \frac{K_p}{\|K_p\|^2} g\circ(\delta \alpha_p) C_{\alpha_p \circ (\delta \alpha_p)}.
\end{equation}
Now
$$
g(\delta\alpha_p(z)) = \frac{1-\bar{p}z}{1-|p|^2\delta - \bar{p}(1-\delta)z};
$$
substituting  the preceding formula for $g\circ(\delta \alpha_p)$ on the right of (\ref{FCRY2}), using $1/\|K_p\|^2 = 1-|p|^2$, and (\ref{phiF}), we obtain (\ref{FCRY}), as desired.
\end{proof}

As a corollary of the proof of the preceding theorem, we obtain a spectral characterization of any normal weighted composition operator on $H^2(\U)$ induced by a selfmap $\phi$ of $\U$ fixing a point in $\U$.

\begin{proposition}\label{SCNIFP} Suppose that $\phi$ has a fixed point  $p\in \U$  and that $W_{\psi, \phi}$ is normal. Then  the spectrum of $W_{\psi, \phi}$ is the closure of  $\{\psi(p)\phi'(p)^n: n = 0, 1, 2, \ldots\}.$
\end{proposition}
\begin{proof}  Because $W_{\psi, \phi}$  is normal and $\phi(p) = p$ for $p\in \U$, we know from the proof of Theorem~\ref{CNWCO} that  $W_{\psi, \phi}$ is unitarily equivalent to $\psi(p)C_{\delta z}$, where $\phi(z) = \alpha_p\circ(\delta \alpha_a)$.  Because $\delta = \phi'(p)$ and the spectrum of $C_{\delta z}$ is the closure of $\{\delta^n: n = 0, 1, 2, \ldots\}$, the unitary equivalence of $W_{\psi, \phi}$ to $\psi(p)C_{\phi'(p)z}$ shows that  the spectrum $W_{\psi, \phi}$ is the closure of  $\{\psi(p)\phi'(p)^n: n\ge 0\}$, as desired. 
\end{proof}

  Continuing with the notation in the proof of the preceding proposition (arising from (\ref{IFPN})), we remark that in case $|\delta| = |\phi'(p)|  < 1$, the operator $W_{\psi, \phi}$ is compact and the spectral characterization  of Proposition~\ref{SCNIFP} follows from \cite[Theorem 1]{GG}.    In case $|\delta| = 1$ and $|\gamma| = |\psi(p)| = 1$, the operator $W_{\psi,\phi}$ is unitary and we obtain the spectral characterization promised by part (a) of Theorem~\ref{SCT}.   
 
   \section{Normality when the inducing function has Denjoy-Wolff point on $\partial \U$}
   
    There are normal weighted composition operators $W_{\psi,\phi}$ where $\phi$ is nonautomorphic and has Denjoy-Wolff point on the unit circle.  In fact, Cowen and Ko \cite{CK} have shown that for every $t >0$ the parabolic map $\phi(z) = T^{-1}(T(z) + t)$, where $T(z) = (1+z)/(1-z)$,  has a companion weight function $\psi$ such that $W_{\psi, \phi}$ is Hermitian.   The weight function $\psi$ in this situation is again linear fractional.   The general form of the weight functions $\psi$  that correspond to  Hermitian weighted composition operators (characterized in \cite{CK}), to unitary weighted composition operators (characterized in Section~\ref{USect} above), and to the normal weighted composition operators (characterized in Section~\ref{SNWCO} above) is easily seen to be given by
\begin{equation}\label{FWF}
    \psi(z) = \rho K_{\sigma(0)}(z)
\end{equation}
    where $\sigma$ is the Cowen axillary function for the linear-fractional inducing map $\phi$ and $\rho$ is a constant.    The following theorem reveals what is required  for normality of a weighted composition operator $W_{\psi, \phi}$, where $\phi$ is linear fractional and $\psi$ has form (\ref{FWF}).   
    
    \begin{proposition}\label{RLFC}  Suppose 
    $$
    \phi(z) = \frac{az+b}{cz+d}
    $$
    is a linear fractional selfmap of $\U$ and $\psi = K_{\sigma(0)}$ where 
    $
    \sigma(z) = (\bar{a}z -\bar{c})/(-\bar{b}z + \bar{d}).$ Then $W_{\psi,\phi}$ is normal
if and only if 
\begin{equation}\label{LFS}
\frac{|d|^2}{|d|^2 - |b|^2 - (\bar{b}a - \bar{d}c)z} C_{\sigma \circ \phi} = \frac{|d|^2}{|d|^2 - |c|^2 - (\bar{b}d - c\bar{a})z} C_{\phi\circ\sigma}.
\end{equation}
\end{proposition}
\begin{proof}  Recall $C_\phi^* = M_g C_\sigma M_h^*$, where $g(z) = 1/(-\bar{b}z + d)$, $h(z) = cz + d$, and $\sigma$ is given in the statement of the proposition.   Note that $K_{\sigma(0)}(z) = \frac{d}{cz + d}$ and $K_{\sigma(0)} h = d$.  Thus,
$$
(W_{\psi,\phi}^* W_{\psi, \phi}) f   =  M_g C_\sigma M_h^* M_{K_{\sigma(0)}}^*  M_{K_{\sigma(0)} }C_\phi \\
                                                 =  \bar{d}\cdot g \cdot K_{\sigma(0)}\circ \sigma \cdot f\circ \phi\circ \sigma,
                                   $$
which simplifies to the right-hand side  of (\ref{LFS}) applied to $f$.   Similarly,
$$
(W_{\psi,\phi} W_{\psi, \phi}^*) f      =  \bar{d} \cdot K_{\sigma(0)} \cdot g\circ\phi \cdot f\circ\sigma\circ \phi.
$$
simplifies to the left-hand side of (\ref{LFS}) applied to $f$.\end{proof}

Given the results of Section~\ref{USect} above and from Section~5 of \cite{CK}, the following result is expected.

\begin{proposition} Suppose that $\phi$ is a linear fractional selfmap of  parabolic type and	$\psi =  K_{\sigma(0)}$, where $\sigma$ is the Cowen auxiliary function for $\phi$; then $W_{\psi,\phi}$ is normal.
\end{proposition}
\begin{proof}   Let $\dw$ be the Denjoy-Wolff point of $\phi$ so that $|\dw| =1$,  $\phi(\dw) = \dw$, and $\phi'(\dw) = 1$.  Just as in the proof of Theorem~\ref{SCT}, by considering conjugation via $C_{\dw z}$, we can, without loss of generality, assume that $\dw =1$.  

Because $\phi$ is parabolic and fixes $1$, it has the from given by (\ref{PF}):
$$
\phi(z) = \frac{(2-t)z + t}{-tz + (2+t)},
$$
where $\Re(t) \ge 0$.   The condition (\ref{LFS}) of Proposition~\ref{RLFC} becomes in this situation
$$
\frac{|2+t|^2}{|2+t|^2 - |t|^2 -4\Re(t)} C_{\sigma\circ \phi} = \frac{|2+t|^2}{|2+t|^2 - |t|^2 -4\Re(t)} C_{\phi\circ \sigma} $$
Because $\phi$ and $\sigma$ have the same fixed point set, they must commute: $\sigma\circ \phi = \phi\circ \sigma$, and the normality of $W_{\psi, \phi}$ follows.
\end{proof}

Suppose that $\phi$ is of hyperbolic type with Denjoy-Wolff point $\dw\in \partial \U$ and $\phi'(\dw) = b < 1$.  Without loss of generality we again assume $\dw =1$ and hence $\phi$ has the form (\ref{HF}): 
$$
\phi(z) =\frac{(1+r-t)z + r+t -1}{(r-t-1)z+1 +r+t},
$$
where $r = 1/b$. Suppose a $\phi$ of this form induces a normal weighted composition operator under the conditions of Proposition~\ref{RLFC}.  Then, applying both sides of (\ref{LFS}) to the constant function $1$ and evaluating the result at $z = 0$, we obtain
$$
\frac{|1+r + t|^2}{|1+r+t|^2 - |r+t-1|^2} =\frac{|1+r + t|^2}{|1+r+t|^2 - |r-t-1|^2}.$$
 It is easy to see that this condition implies $\Re(t) = 0$ so that $\phi$ is an automorphism.  Thus, no hyperbolic non-automorphic linear fractional map (with $\dw \in \partial \U$) can induce a normal weighted composition operator under the conditions of Proposition~\ref{RLFC}.   Some evidence that this is true in general may be found in \cite{BLNS}, whose results show that 
 \begin{quotation}
 $W_{\psi, \phi}$ cannot even be essentially normal if $\psi$ is, say, $C^1$ on the closure of $\U$ and  $\phi$ is  linear-fractional nonautomorphism with Denjoy-Wolff point $\dw\in \partial \U$ and $\phi'(\dw) < 1$.
 \end{quotation}
  Suppose $\psi$ is $C^1$ on the closure of $\U$ and $\phi$ is a linear-fractional nonautomorphism having Denjoy-Wolff point $\dw \in \partial \U$. Then Lemma 3.3 of \cite{BLNS} shows that $W_{\psi, \phi}$ is equivalent to $\psi(\dw)C_\phi$ modulo the compact operators; moreover, if $\phi'(\dw) < 1$, then Theorem 5.2 of \cite{BLNS} shows that $\psi(\dw)C_\phi$ is not essentially normal, so that in this situation $W_{\psi, \phi}$ is not essentially normal. 
  


\begin{thebibliography}{99}
              
 \bibitem{BLNS} P.~S~Bourdon, D.~Levi, S.~K.~Narayan, and J.~H.~Shapiro,  Which linear-fractional composition operators are essentially normal?, {\it  J. Math. Anal. App.}\  280 (2003) 30--53.

\bibitem{BrdS} P.~S.~Bourdon, Spectra of some composition operators and associated weighted composition operators, preprint.  


\bibitem{Cow2}  C.~C.~Cowen, Composition operators on $\htwo$, {\it  J.\ Operator Th.}\ 9 (1983),  77--106.

\bibitem{Cow} C.~C.~Cowen, Linear fractional composition
             operators on $\htwo$, {\it  Integral Eqns.\ Op.\ Th.}\  11 (1988),
              151--160.
              

\bibitem{CK} C.~C.~Cowen and E.~Ko, Hermitian weighted composition operators on $H^2$, {\it Trans.\ Amer.\ Math.\ Soc.}, to appear.


\bibitem{CMB} C.~C.~Cowen and B.~D.~MacCluer, {\em Composition
Operators on Spaces of Analytic Functions\/}, CRC Press 1995.



\bibitem{GG} G.~Gunatillake, Spectrum of a compact weighted composition operator,  {\it Proc.\ Amer.\ Math.\ Soc.}\  135  (2007), 461--467.


\bibitem{Lit} J.~E.~Littlewood,  On inequalities in the
theory of functions, {\it  Proc.\ London Math.\ Soc.}\ 23 (1925),
481--519.

\end{thebibliography}
\end{document}